\documentclass[aps,twocolumn,superscriptaddress]{revtex4}
\usepackage{amsmath,amssymb}
\usepackage{graphics,graphicx}
\usepackage{dcolumn,bm}

\newcommand{\imsc}{\affiliation{The Institute of Mathematical
Sciences,
CIT Campus, Taramani, Chennai 600 113, India.}}
\newcommand{\tifr}{\affiliation{Department of Astronomy and
Astrophysics,
Tata Institute of Fundamental Research, Homi Bhabha Road, Mumbai 400
005,
India.}}

\begin{document}
\title{In Square Circle: Geometric Knowledge of the Indus
Civilization}
%\title{Circular geometry and space-filling in Indus Valley
%civilization artfifacts}

\author{Sitabhra Sinha}%
\email{sitabhra@imsc.res.in}
\imsc
\author{Nisha Yadav}%
\tifr
\author{Mayank N. Vahia}%
\tifr

%\date{\today}

\begin{abstract}
The earliest origins of mathematics in the Indian subcontinent is
generally dated around 800-500 BCE when the {\em Sulbasutras}
are thought to have been written. In this article we suggest that
mathematical thinking in South Asia, in particular, geometry, may have 
had an even earlier beginning - in the third millenium BCE. 
We base our hypothesis on the analysis of
design patterns, such as complex space-filling tiling,
seen on artifacts of the Indus Valley
Civilization (also referred to as the 
Mature Harappan Civilization, 2500-1900 BCE) 
which speaks of a deep understanding of
sophisticated geometric principles. 
\end{abstract}

\maketitle
\section{Introduction.}
The geometric principles expounded in the {\em Sulbasutras} (800-500
BCE) have often been considered to mark the beginning of mathematics
in the Indian subcontinent~\cite{Seidenberg75,Staal08}. 
This collection of {\em sutras} codify directions for constructing
sacrificial fires, including rules for the complex
configuration of ritual altars. The bird-shaped {\em Agnicayan} altar, 
consisting of five layers of two hundred bricks each, with the bricks
being of square,
rectangular or triangular shapes of various sizes, is considered by
F. Staal to signal the beginning of geometry proper in South
Asia~\cite{Staal08}. It has been dated by him to about 1000 BCE as
some of the mantras that are concerned with the consecration of bricks
occur in the earliest Yajurveda Samhita, the Maitrayani.
The absence of any recorded tradition of geometric knowledge predating
these sutras have led some scholars to suggest a West Asian origin
for the onset of mathematical thinking in India. However, the
discovery of the archaeological remnants of the Indus Valley
civilization in parts of Pakistan and northwestern India over the
course of last century has revealed a culture having a sophisticated 
understanding of geometry which predated the {\em Sulbasutras} by more
than a thousand years. It is difficult to ascertain whether there was
any continuity between the geometry practised by the
Indus civilization and that used by the later Vedic culture; however,
it is not impossible that some of the earlier knowledge persisted
among the local population and influenced the {\em sulbakaras}
(authors of the {\em Sulbasutras}) of the
first millennium BCE.

\section{Indus geometry: the archaeological evidence}
The Indus Valley civilization, also referred to as the Mature Harappan
civilization (2500-1900 BCE), covered approximately a million square
kilometres -- geographically spread over what are now Pakistan and 
northwestern India \cite{Possehl02}. It was approximately
contemporaneous with Old Kingdom Egypt, as well as, the Early Dynastic
Period of the Sumerian Civilization and the Akkadian Empire in
Mesopotamia. The Indus culture was characterized by extensive urbanization 
with large planned cities, as seen from the ruins of Harappa
and Mohenjo-daro (among others). There is evidence for craft
specialization and long-distance trade with Mesopotamia and Central
Asia.

The well-laid out street plans of the Indus cities and their 
orientation along the cardinal directions have been long been taken as
evidence that the Indus people had at least a working knowledge of
geometry~\cite{Amma79,Parpola94}. Earlier studies have suggested that
not only did these people have a practical grasp of
mensuration, but that they also had an understanding of the basic
principles of geometry~\cite{Kulkarni78}. The discovery of scales and
instruments for measuring length in different Indus sites indicate
that the culture knew how to make accurate spatial measurements
\cite{Vij84,Balasubramaniam08}. 
For example, an
ivory scale discovered at Lothal (in the western coast of India) 
has 27 uniformly spaced lines over 46 mm, indicating an unit of length
corresponding to 1.70 mm~\cite{Rao85}. 
The sophistication of the metrology 
practised by the Indus people is attested by the sets of regularly shaped
artifacts of various sizes that have been identified as constituting a system
of standardized weights.
Indeed, there is a surprising degree of uniformity in the measurement units
used at the widely dispersed
centers of the civilization, indicating an attention towards
achieving a standard system of units for measurement.

However, most of the literature available up till now on the subject of 
Indus geometry have been primarily concerned with patterns occurring
at the macro-scale (such as building plans, street alignments, etc.).
The smaller-scale geometric patterns that are often observed on seals or on the
surface of pottery vessels have not been analysed in great detail.
We believe that such designs provide evidence for a much more
advanced understanding of geometry on the part of the Indus people
than have hitherto been suggested in the literature.

\begin{figure}
\includegraphics[width=0.7\linewidth]{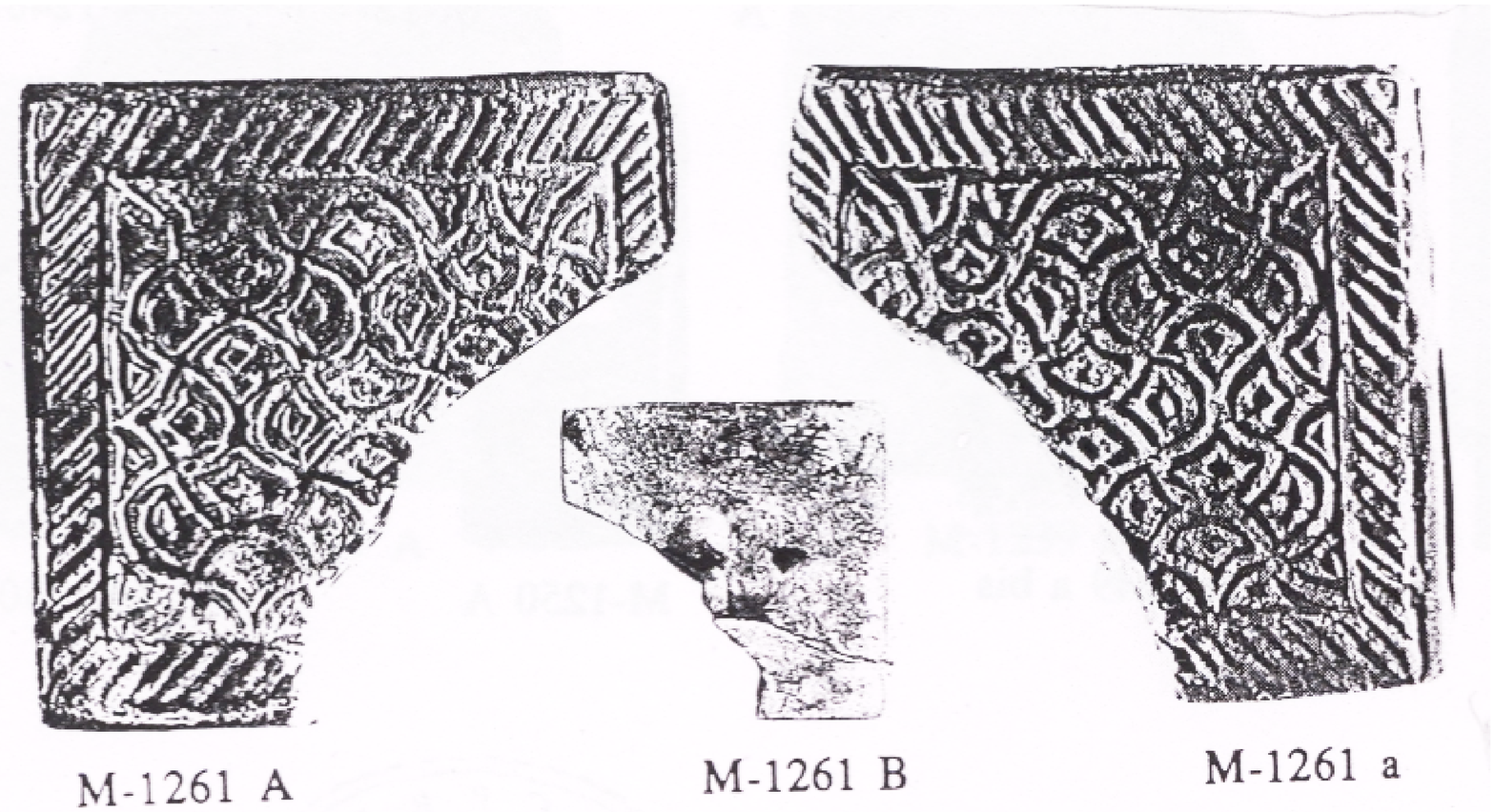}
\includegraphics[width=0.29\linewidth]{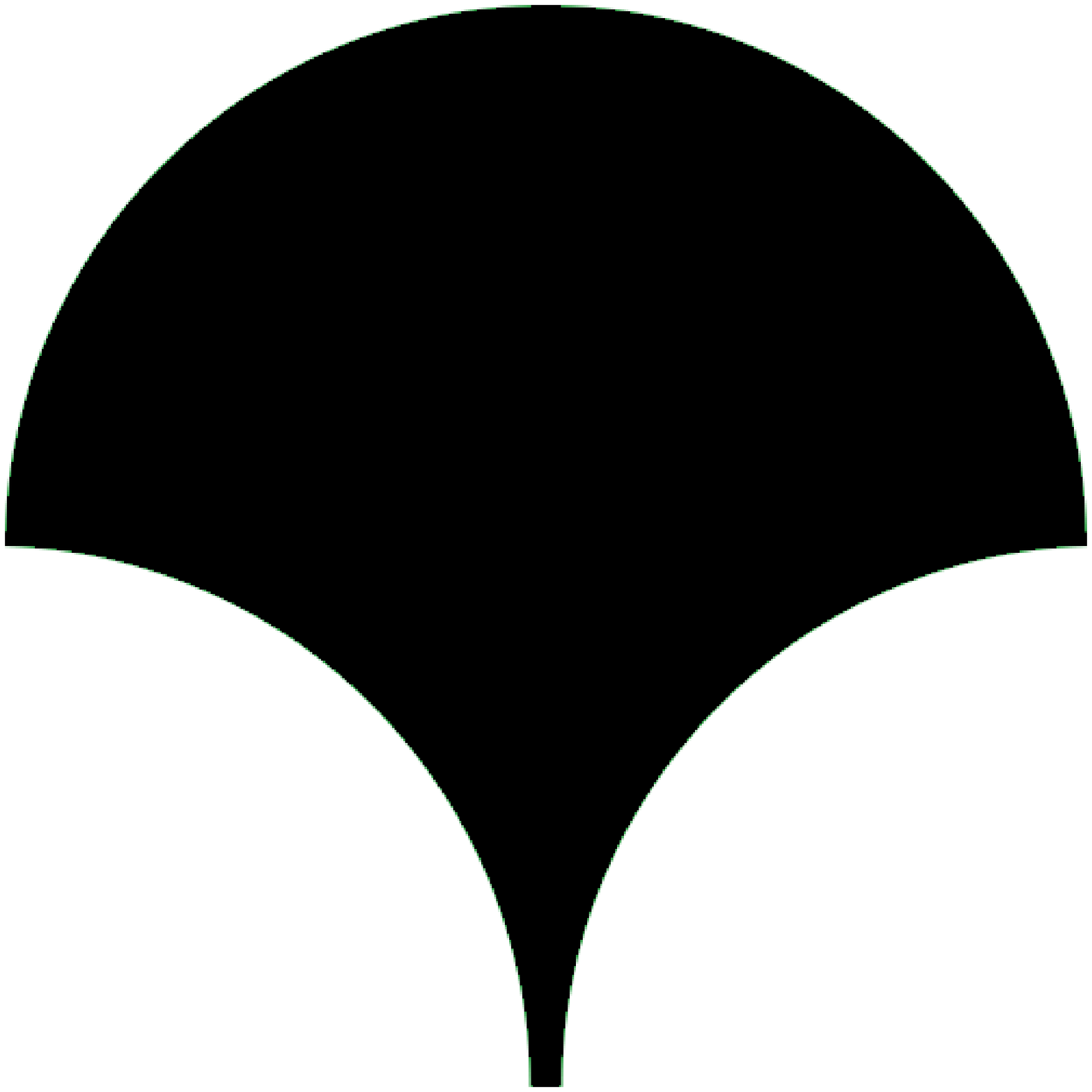}
\caption{\label{fig:M_1261}
Space-filling pattern in a seal excavated from Mohenjo-daro, 
DK Area (No. M-1261 in p. 158 of Ref.~\cite{Shah91}) showing the 
seal (A), its impression (a) and the boss side (H). On the right, we
show the basic motif having the shape of a Japanese fan.
[Figure of the seal is reproduced with permission from A. Parpola.]
}
\end{figure}

\begin{figure}
\includegraphics[width=0.95\linewidth]{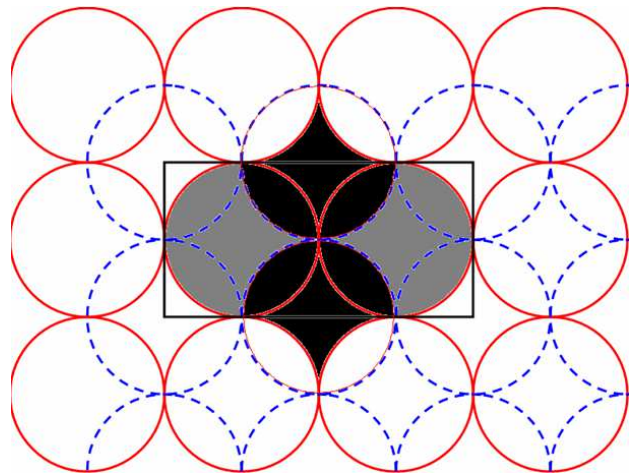}
\caption{\label{fig:space_filling_circles}
Close-packing of circles in a four-fold symmetric
arrangement having two layers. The circles in the different layers are
distinguished by broken and continuous curves, respectively.
The two layers are displaced with
respect to each other by half a lattice spacing in both the horizontal
and vertical directions.
The marked out rectangle shows the area
which is transformed into the shape indicate by shaded areas.
}
%\vskip 1in
\end{figure}

\begin{figure}
\includegraphics[width=0.95\linewidth]{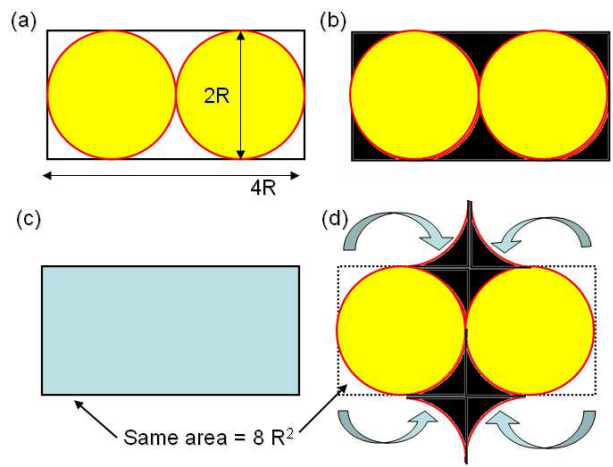}
\caption{\label{fig:transformation}
The set of transformations that produce the repeating motif for the
space-filling tiling of Fig.~\ref{fig:M_1261}. Two non-overlapping
circles, each of radius $R$ are inscribed within a rectangle of length
$4R$ and width $2R$ (a). Of the sections of the rectangle not belonging
to the circles (shaded black in (b)), four are cut and re-pasted
so that they fall outside the original rectangle (c), creating the
repeating motif. 
Note that, its area ($= 8 R^2$) is the same as that of the original rectangle 
(d). 
}
%\vskip 1in
\end{figure}

As an example, we direct the reader's attention to the occurrence of a
space-filling pattern on a seal, where a basic motif shaped like a
Japanese fan has been used in different orientations to tile an
approximately square domain (Fig.~\ref{fig:M_1261}). 
The apparent simplicity of the design
belies the geometric complexity underlying it. The actual repeating pattern
is a fairly complicated shape that is composed of four of the basic
fan-shaped motifs, each being rotated by 90 degrees (see the shaded
object in Figure~\ref{fig:space_filling_circles}). The
transformations that are required to construct this pattern from a
rectangle whose length is twice its breadth is shown in
Figure~\ref{fig:transformation}. Note that, despite the complex shape
of the repeating pattern, its area is fairly simple to calculate.
By construction, this area is equal to that of the rectangle. The
construction also
enables us to calculate the area of the fan-shaped basic motif as
being $2R^2$ where $R$ is the length from its sharply pointed tip to
the midpoint of its convex edge.

The set of operations that are required to generate the shapes
described above suggest that the people of the Indus civilization 
were reasonably
acquainted with the geometry of circular shapes and techniques of
approximating their areas, especially as we know of several Indus
artifacts which exhibit other designs that follow from these
principles. In fact, we note that the Indus civilization paid 
special attention to the circle and its variants in the geometric designs 
that they made on various artifacts.
It reminds us of the special place that the circular shape had
in Greek geometry that is aptly summarised by the statement of 
the 5th century neo-Platonist Greek philosopher Proclus in his
overview of Greek geometry: ``The first and simplest
and most perfect of the figures is the circle'' \cite{Proclus}.

The reason for the primacy of the circle in Indus geometry is probably not 
hard to understand if we focus on the technical means necessary to
generate ideal geometric shapes. Unlike rectilinear shapes which are
very difficult to draw exactly without the help of instruments that
aid in drawing lines
at right angles to each other, an almost perfect circle can be
drawn by using a rope and a stick on soft clay. Indeed there is
evidence for the use of compass-like devices for drawing circles in the
Indus Valley civilization. E. J. H. Mackay, who excavated Mohenjo-daro
between 1927-1932 expressed surprise on finding that ``an instrument was
actually used for this purpose [drawing circles] in the Indus Valley
as early as 2500 BC"~\cite{Mackay38}. In this context, it may be noted
that later excavations at Lothal have unearthed thick, ring-like shell
objects with four slits each in two margins that could have been used
to measure angles on plane surfaces~\cite{Rao85}. 

\begin{figure}[htb]
\includegraphics[width=0.95\linewidth]{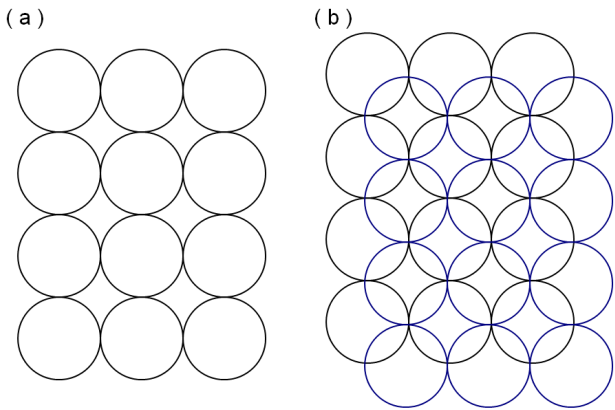}
\caption{\label{fig:circle_lattice}
(a) A periodic lattice consisting of close-packed circles in a
four-fold symmetric
arrangement. (b) The intersecting circle lattice generated by superposing
two lattices of circles displaced with respect to each other by half a lattice
spacing in the vertical and horizontal directions.}
%\vskip 1in
\end{figure}

Let us return to the issue of how the sophisticated understanding of
tiling the plane with shapes having circular edges may have emerged 
in the Indus
civilization. A possible clue can be found from the construction given
in Fig.~\ref{fig:space_filling_circles}, viz., the close-packing of
circles on a plane in a four-fold symmetric arrangement. Such lattice patterns
with the spatial period equal to the diameter of a circle, are
easy to generate, for instance, by pressing an object with a circular 
outline on a wet clay surface repeatedly along rows and columns
(Fig.~\ref{fig:circle_lattice}~(a)). 
Indeed, certain Indus artifacts bear the impression of circular objects pressed
on them to create a design. By exploring the variations of such a basic
pattern, a variety of increasingly complicated design motifs can be generated
which culminate in the complicated space-filling tiles shown in 
Figure~\ref{fig:M_1261}.

\begin{figure}
\includegraphics[width=0.45\linewidth]{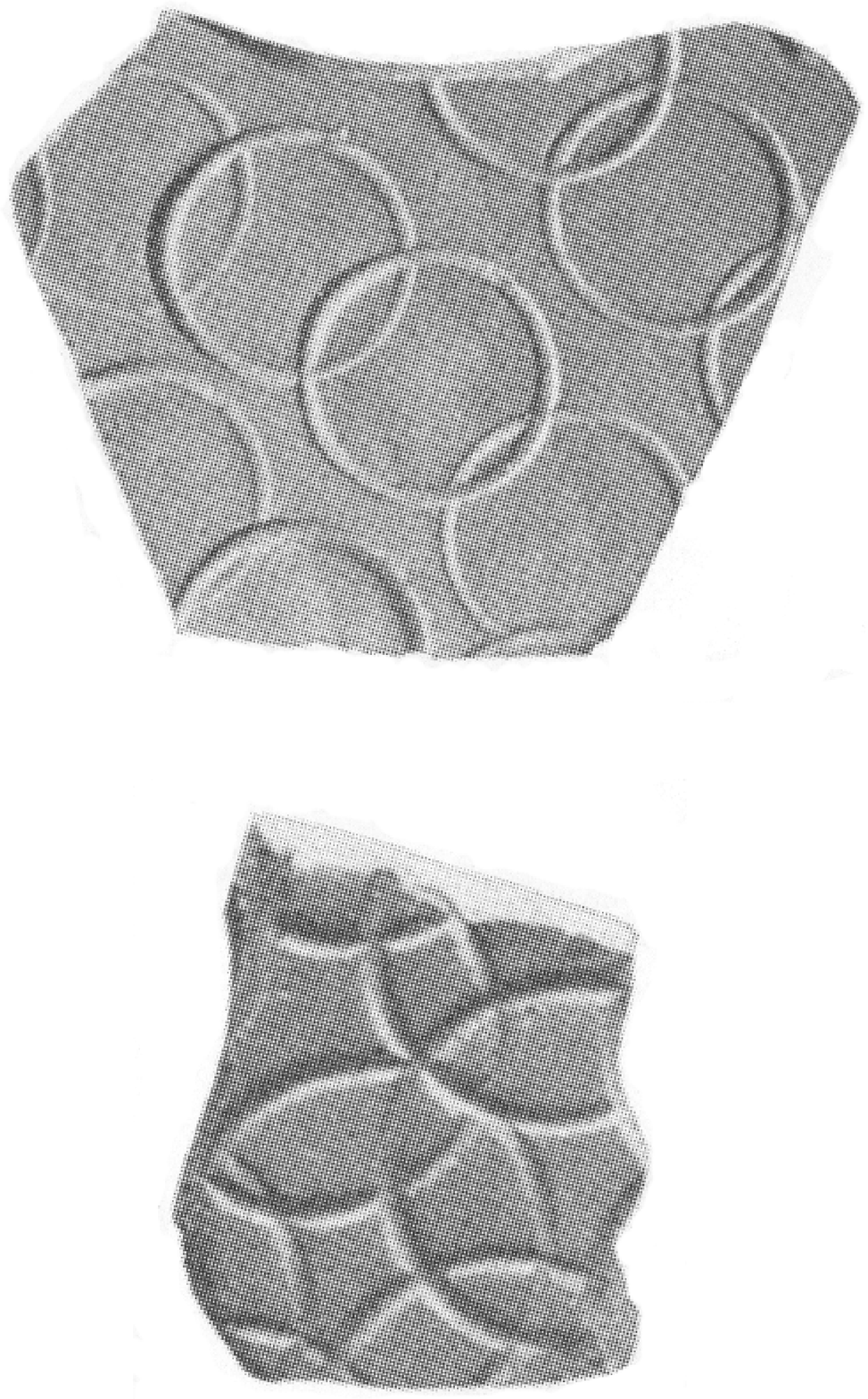}
\includegraphics[width=0.45\linewidth]{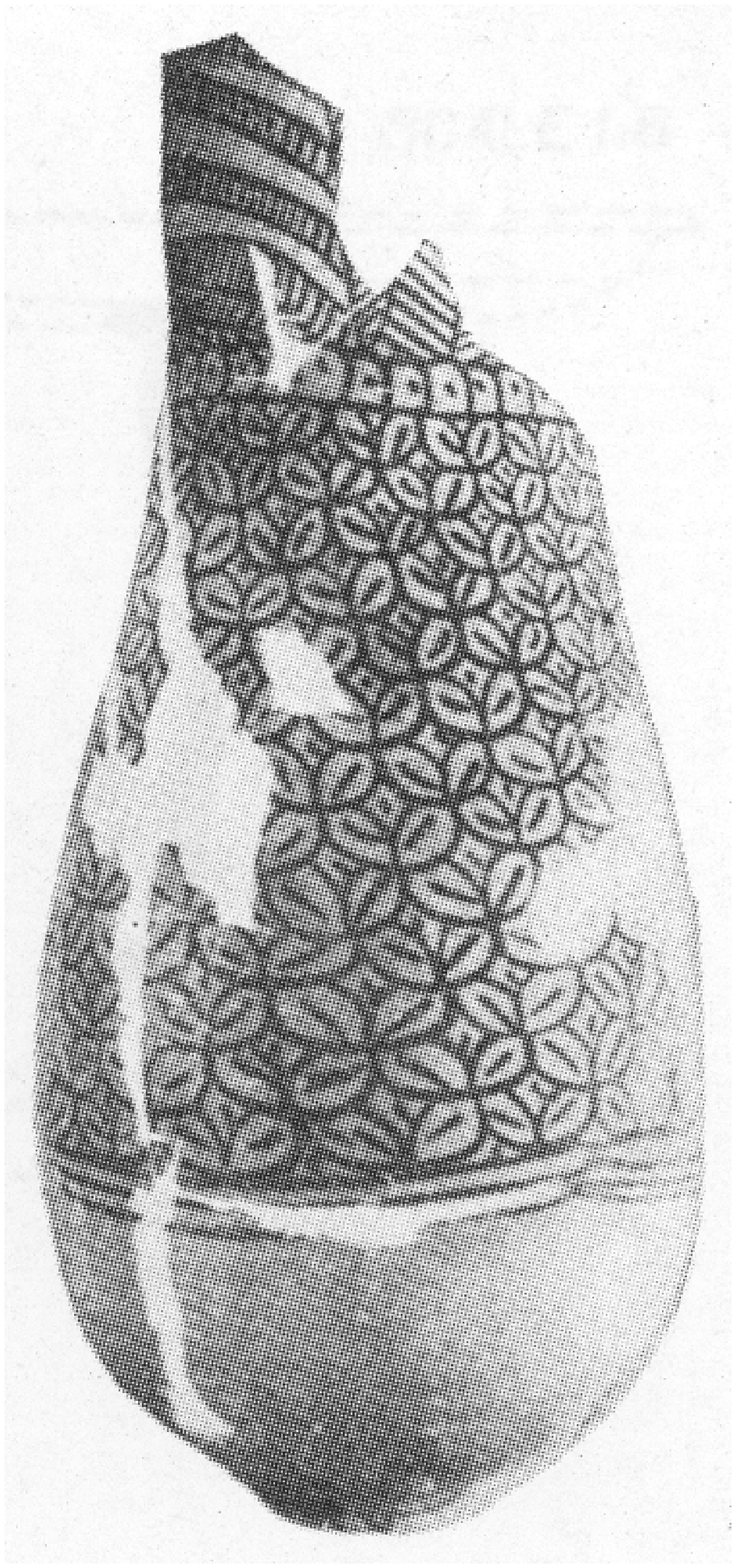}
\caption{\label{fig:intersecting_circles}
(Left) The motif of intersecting circles on incised pottery from lower levels
of
Mohenjo-daro, DK Area, G section (Nos. 24 [top] and 21 [bottom]
in Plate LXVII of Ref.~\cite{Mackay38}). (Right) A decorated
intersecting circle lattice used as a space filling pattern
on the surface of jar from Mohenjo-daro, DK Area, G section (no. 5
in
Plate LIV of Ref.~\cite{Mackay38}).}
%\vskip 1in
\end{figure}

 One of the simplest variations of the basic circular pattern is to have two
layers of such close-packed lattices of circles, with one layer displaced by
a length equal to the radius of a circle (i.e., half the lattice spacing) along both the horizontal
and vertical directions (Fig.~\ref{fig:circle_lattice}~(b)). The resulting    
intersecting circle lattice is a frequently occurring design motif in
Indus artifacts. See Figure~\ref{fig:intersecting_circles} 
for examples of this pattern,
which occasionally occurs in a ``decorated" variation (with points or lines used
to fill the spaces within the circular shapes). 
In this context, one may mention that one of the signs that
occur in frequently in the Indus seal inscriptions is in
the form of two overlapping ellipses.
Parpola has suggested as association of this sign with the Pleiades star
system~\cite{Parpola94}. It is a matter of conjecture whether the
design of intersecting circles found in so many Indus artifacts has any
astral significance.

\begin{figure}[htb]
\begin{center}
\includegraphics[width=0.66\linewidth]{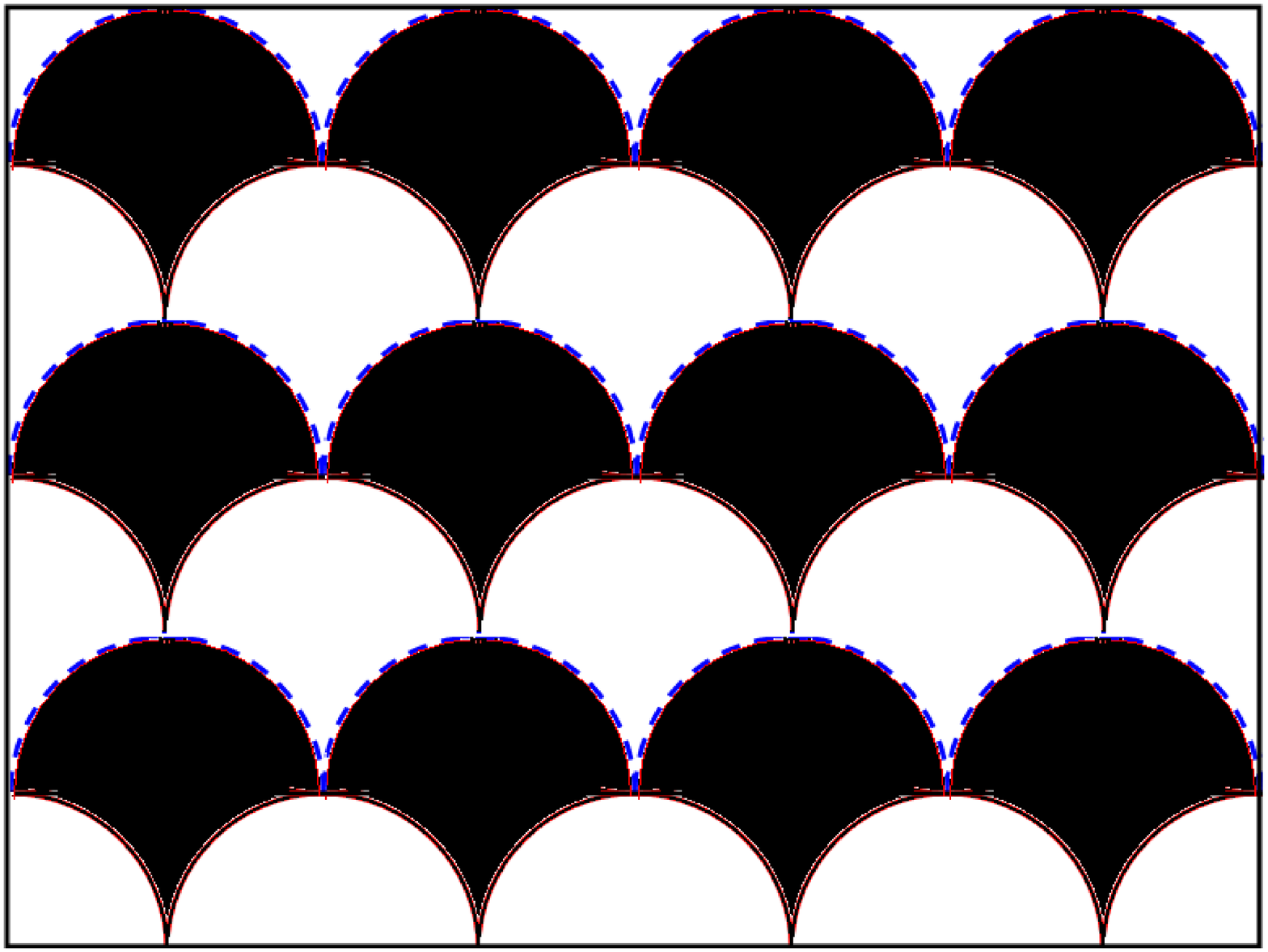}
\includegraphics[width=0.275\linewidth]{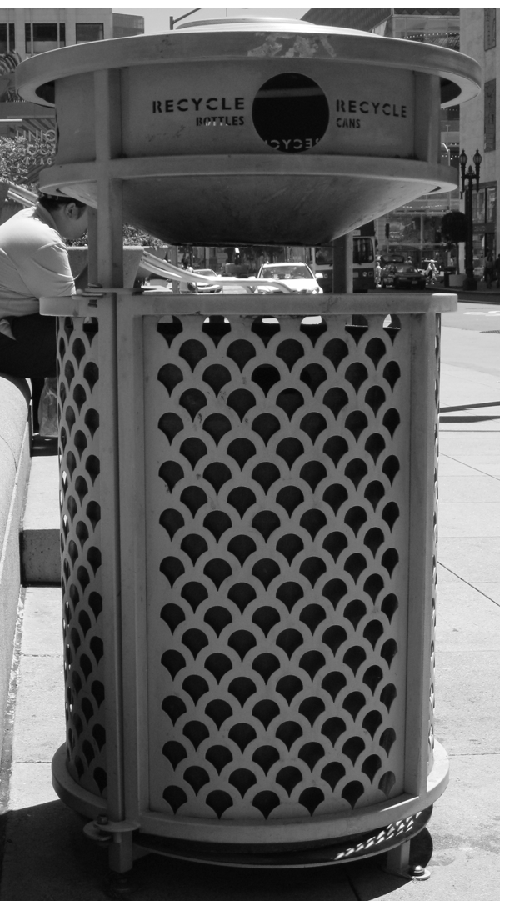}
\end{center}
\caption{\label{fig:fish_scale}
(Left) Imbricate pattern of overlapping circles used for tiling a plane. The basic
fan-shaped motif is identical to the one show in 
Fig.~\ref{fig:M_1261}. (Right) The same pattern as it appears in a present-day 
urban setting (Photo: Nivedita Chatterjee).}
%\vskip 1in
\end{figure}
\begin{figure}
\includegraphics[width=0.95\linewidth]{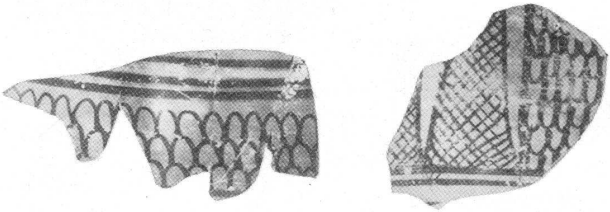}
\caption{\label{fig:fish_scale_actual}
Imbricate pattern of overlapping circles seen on fragments of painted pottery 
excavated from
Mohenjo-daro,
DK Area, G section (Nos. 33-34 in Plate LXX of Ref.~\cite{Mackay38}).}
%\vskip 1in
\end{figure}
The next variation we discuss is the imbricate pattern having regularly arranged overlapping
edges resembling fish scales. This pattern has also frequently appeared in artifacts of 
other cultures, including on the surface of  
West Asian Intercultural Style vessels~\cite{Kohl78}. 
Figure~\ref{fig:fish_scale} shows
how the design is used for tiling the plane by having several rows of overlapping circles,
each partially obscuring the layer beneath it. It is fairly easy to see that
this is a variation of the intersecting circle lattice. One can envisage the design as being
made up of multiple rows of circular tiles stacked one on top of the other with successive rows
displaced by half a lattice spacing vertically and horizontally.
The pattern can be seen frequently on the surface of Indus painted
pottery,
both in its original and ``decorated" forms 
(Fig.~\ref{fig:fish_scale_actual})~\cite{Starr41}.

\begin{figure}
\includegraphics[width=0.95\linewidth]{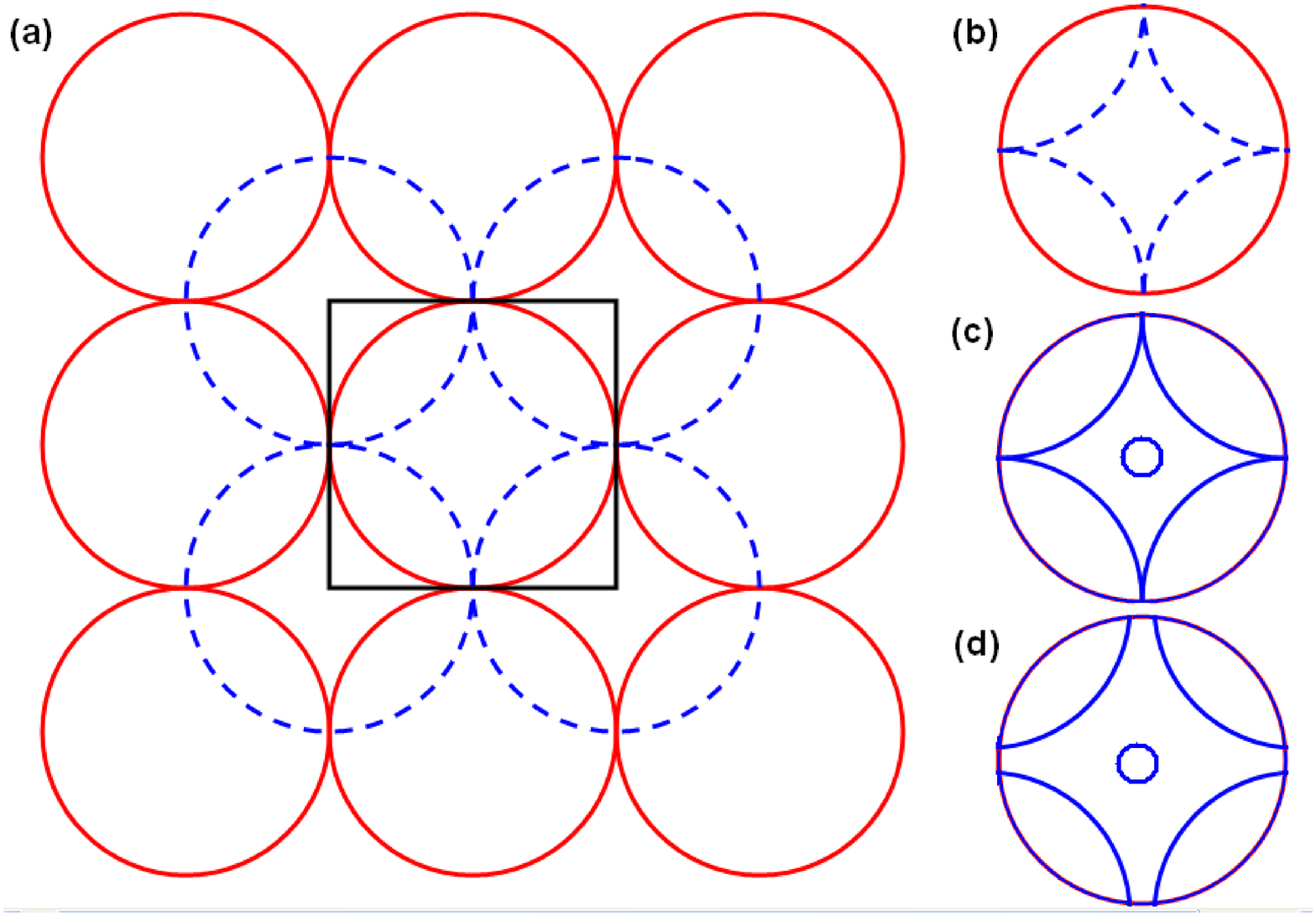}
\caption{\label{fig:quincross}
The intersecting circle lattice (a) can be used to generate a star-like motif (b)
with decorated variants (c) and (d).}
%\vskip 1in
\end{figure}

\begin{figure}[tbp]
\includegraphics[width=0.25\linewidth]{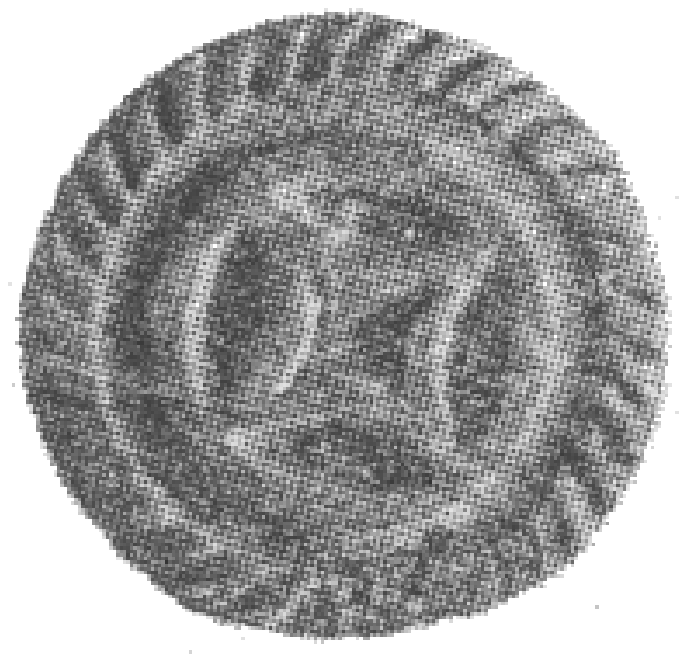}
\includegraphics[width=0.25\linewidth]{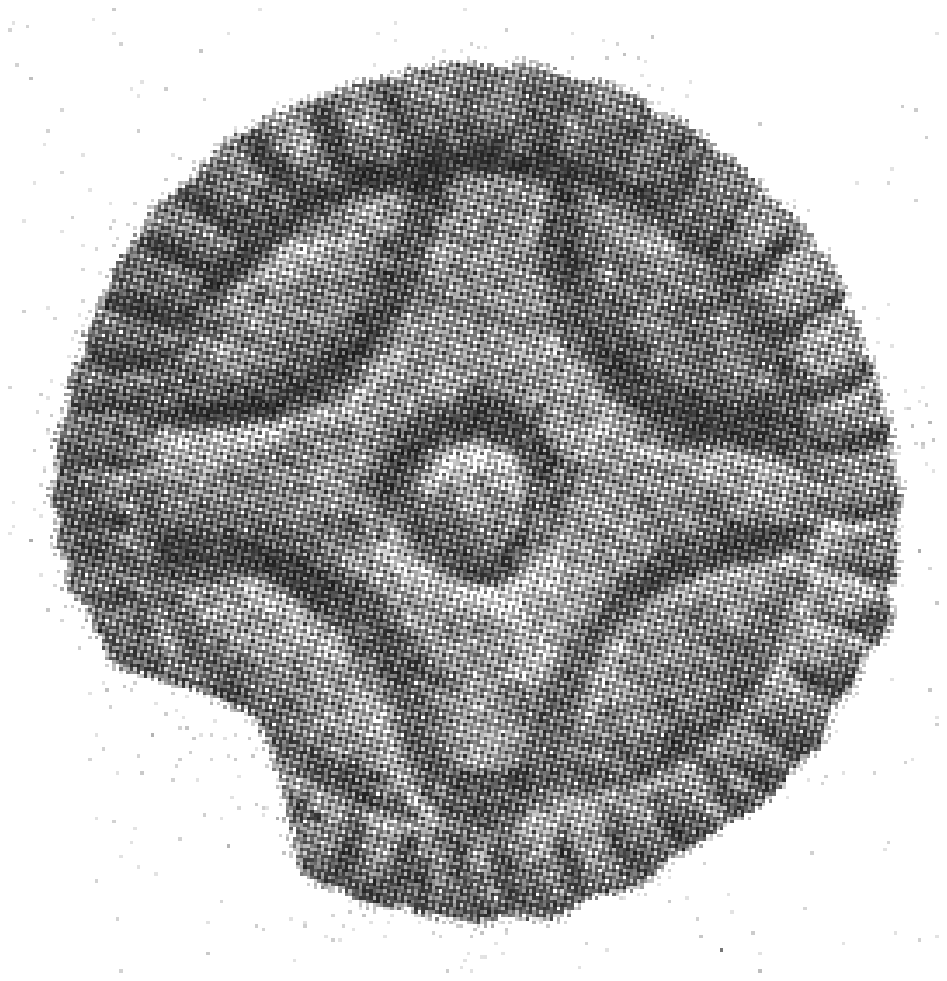}
\includegraphics[width=0.45\linewidth]{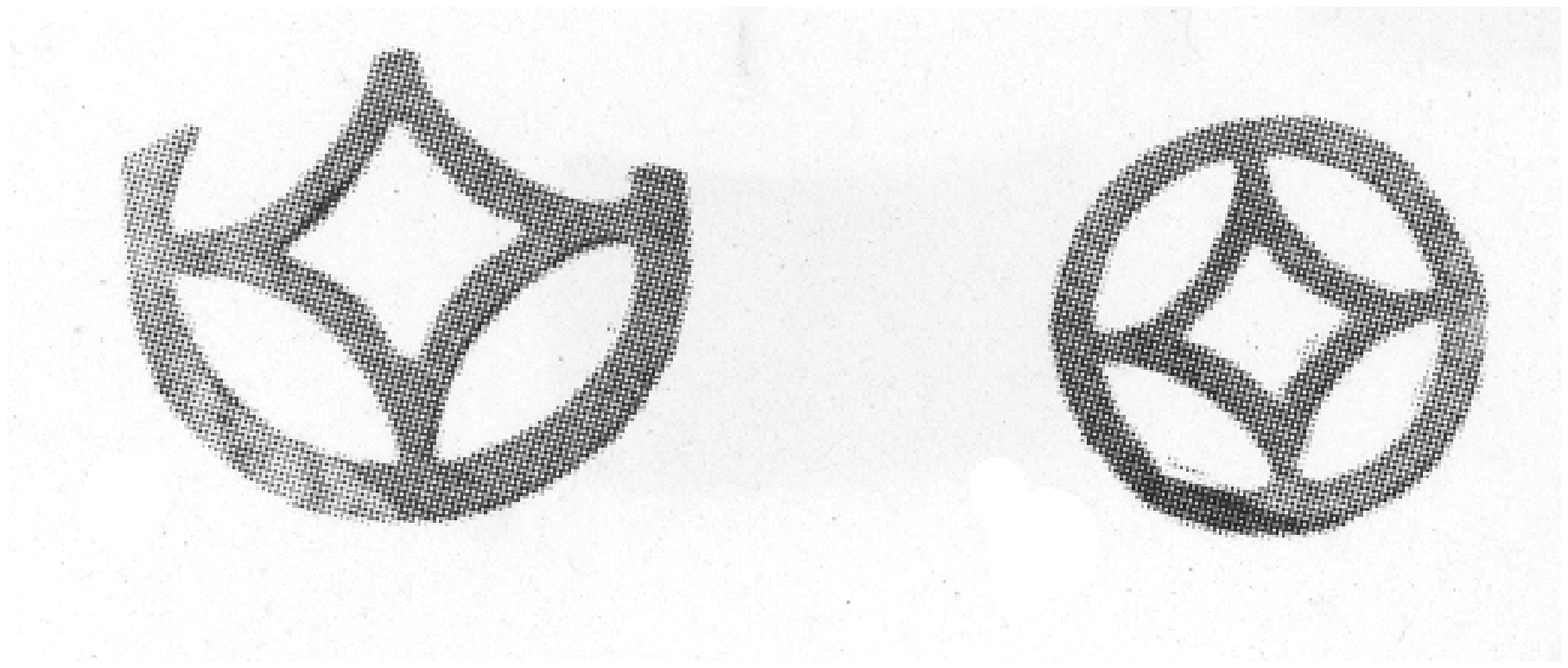}
\caption{\label{fig:quincross_actual}
The star-like geometrical motif in
ornament-like objects excavated from Harappa (left) and in a
seal impression (center) and ivory objects (right) found in
Mohenjo-daro,
DK Area,
G section. (No. 8 in Plate CXXXIX of Ref.~\cite{Vats40}, No. 12
in
Plate C [center] and Nos. 9-10 in Plate CXLI [right]
of Ref.~\cite{Mackay38}).}
%\vskip 1in
\end{figure}

Yet another possible variation obtained from the intersecting circle lattice
is shown in Figure~\ref{fig:quincross}.  If a square, the length of whose sides
equal the diameter of each circle, is placed on the lattice so as to completely
contain a circle inside it, one obtains a star-shaped motif by focusing exclusively within
the region enclosed by the square (Fig.~\ref{fig:quincross}~(b)). A variant
of this four-pointed star is obtained by relaxing the close-packed nature of
the lattice so that the circles no longer touch each other.
The resulting `cross'-shaped motif with a small circle added at the center 
(Fig.~\ref{fig:quincross}~(d)), sometimes referred to as the quincross
sign,
is found to occur very frequently in pre-Columbian iconography,
especially in Mayan (`Kan-cross') and Zapotec (`Glyph E')
designs~\cite{Sugiyama05}. The meaning attached to this motif
has been interpreted differently by various archaeologists,
and in the context of the Maya, it has been claimed to be a
representation of Venus. A possible astral significance of this sign in the
context of Indus civilization is probably not far-fetched
given its star-like shape. The motif, in both the exact four pointed star-like shape
and in its cross-like variation, is found in many artifacts at both Harappa and
Mohenjo-daro (Fig.~\ref{fig:quincross_actual}).

The intersecting circle lattice is thus a fairly general mechanism for
producing a variety of geometric patterns, ranging from quite simple
ones that have been seen to occur across many different cultures to
complicated tiling shapes that suggest a deeper understanding of the
underlying geometry. 
This is especially the 
case for the complicated space-filling tile made from the four fan-shaped
motifs shown in Fig.~\ref{fig:transformation}. 
It is clear from the construction procedure that the geometrical knowledge
underlying the origin of this design is non-trivial and implies a sophisticated
understanding of the principles of circular geometry. 

\section{Conclusion}
To summarize, in this article we have argued that the origin of
mathematics, and geometry in particular, in the Indian subcontinent 
may actually date back to the third millennium BCE and the Indus Valley
civilization, rather than
beginning with the {\em Sulbasutras} of the first millennium BCE as is
conventionally thought. Although the well-planned cities, standardized
system of measurements and evidence of long-range trade contacts with
Mesopotamia and Central Asia attest to the high technological
sophistication of this civilization, the unavailability of written
records has up to now prevented us from acquiring a detailed understanding 
of the level of mathematical knowledge attained by the Indus people.
By focusing on the geometry of design motifs observed commonly in the
artifacts excavated from various sites belonging to this culture, we
have shown that they suggest a deep understanding of the properties of
circular shapes. In particular, designs which exhibit space-filling tiling 
with complicated shapes imply that the Indus culture may have
been adept at accurately estimating the area of shapes enclosed by
circular arcs. This speaks of a fairly sophisticated knowledge of
geometry that may have persisted in the local population long after the
decline of the urban centers associated with the civilization and it
may well have influenced in some part the mathematical thinking of the 
later Vedic culture of the first millennium BCE.

\noindent
{\bf Acknowledgements:}

\noindent
We would like to thank P. P. Divakaran and F. Staal for several 
insightful comments and suggestions on an early version of the manuscript.
S.S. would like to gratefully acknowledge Dr P. Pandian and the staff
of IMSc Library for their help in obtaining several rare
volumes describing the Indus civilization artifacts.

\end{document}